\documentclass{amsart}
\font\fiverm=cmr5
\font\ehsc=cmcsc10 scaled 850

\let\sse=\subseteq
\let\noi=\noindent
\let\veps=\varepsilon
\let\limply=\Longrightarrow

\def\0{\{0\}}
\def\span{{\kern.5pt{\rm span}\kern1pt}}

\def\conv{{\;\longrightarrow\;}}
\def\wconv{{{\buildrel_{\scriptstyle w}\over\conv}}}
\def\sconv{{{\buildrel_{\scriptstyle s}\over\conv}}}
\def\uconv{{{\buildrel_{\scriptstyle u}\over\conv}}}
\def\sslash{\hbox{{\fiverm/}}}
\def\notconv{{{\conv\kern-13pt\slash}\kern9pt}}
\def\notuconv{{{\uconv\kern-13pt\sslash}\kern9pt}}
\def\notsconv{{{\sconv\kern-13pt\sslash}\kern9pt}}
\def\notwconv{{{\wconv\kern-13pt\sslash}\kern9pt}}
\def\query{{{\buildrel_{^{\scriptstyle ?}}\over\limply}}}

\def\B{{\mathcal B}}

\def\N{{\mathcal N}}
\def\Oe{{\mathcal O}}
\def\R{{\mathcal R}}
\def\X{{\mathcal X}}
\def\BX{{\B[\X]}}

\def\CC{{\mathbb C\kern.5pt}}
\def\FF{{\mathbb F\kern.5pt}}
\def\RR{{\mathbb R\kern.5pt}}

\newsymbol\varnothing 203F
\let\void=\varnothing

\def\matrix#1{\null\,\vcenter{
              \normalbaselines\mathsurround=0pt\ialign{
              \hfil $##$
              \hfil && \quad
              \hfil $##$
              \hfil \crcr
              \mathstrut \crcr
              \noalign{\kern-\baselineskip}#1 \crcr
              \mathstrut \crcr
              \noalign{\kern-\baselineskip} \crcr }}\,}

\begin{document}

\vglue-65pt\noi
\hfill{\it Rendiconti del Circolo Matematico di Palermo}\/,
{\bf 73}(2) (2024) 663--673

\vglue15pt
\title[Weak Supercyclicity and Quasistability]
{Weak l-Sequential Supercyclicity and Weak Quasistability}
\author{C.S. Kubrusly}
\address{Catholic University of Rio de Janeiro, Brazil}
\email{carlos@ele.puc-rio.br}
\author{B.P. Duggal}
\address{Faculty of Sciences and Mathematics, University of Ni\v s, Serbia}
\email{bpduggal@yahoo.co.uk}
\renewcommand{\keywordsname}{Keywords}
\date{June 13, 2023}

\begin{abstract}
It is known that supercyclicity implies strong stability$.$ It is not known
whether weak l-sequential supercyclicity implies weak stability$.$ In this
paper we prove that {\it weak l-sequential supercyclicity implies weak
quasistability}$.$ Corollaries concerning the characterisation of (i) weakly
l-sequentially supercyclic vectors that are not (strongly) supercyclic, and
(ii) weakly l-sequentially supercyclic isometries, are also proved.
\end{abstract}

\maketitle

\vskip-20pt\noi
\section{Introduction}

\vskip6pt
Supercyclicity implies strong stability$.$ That is, if the projective orbit
${\Oe_T([y])}$ of a power bounded linear operator $T$ on a normed space $\X$
is dense in $\X$ for some vector $y$ in $\X$, then ${T^nx\to0}$ for every
${x\in\X}$ (both density and convergence are in the norm topology)$.$ In other
words, ${\Oe_T([y])^-\!=\X}$ implies ${T^n\!\sconv O}.$ This is an important
result from \cite[Theorem 2.2]{AB}$.$ It has been asked in \cite{KD1} whether
such an implication survives from norm topology to weak topology$.$ In
particular, whether a stronger version of weak supercyclicity (which is weaker
than supercyclicity in the norm topology, but stronger than supercyclicity in
the weak topology) implies weak stability (i.e., implies that ${T^nx\wconv0}$
for every ${x\in\X}).$ Such a stronger version of weak supercyclicity is
called weak l-sequential supercyclicity, which means that the weak limit set
$\Oe_T([y])^{-wl}$ of the projective orbit (i.e., the set of all weak limits
of weakly convergent sequences in the projective orbit) coincides with $\X.$
So it was asked whether ${\Oe_T([y])^{-wl}\!=\X}$ implies ${T^n\!\wconv O}$,
where weak stability is defined by the expression: ${\lim_n|f(T^nx)|=0}$ for
every ${x\in\X}$, for every ${f\in\X^*}.$ A weaker form of weak stability,
called weak quasistability, is defined as: ${\liminf_n|f(T^nx)|=0}$ for every
${x\in\X}$, for every ${f\in\X^*}.$ In \cite[Theorem 6.2]{KD1} the authors
gave a sufficient condition for a weakly l-sequentially supercyclic power
bounded operator to be weakly stable$.$ Here we improve that result by showing
that, under the power boundedness assumption,
$$
\hbox{\it weak l-sequential supercyclicity implies weak quasistability}.
$$
Notation and terminology used above will be defined in the next section$.$ The
paper is organized into five sections$.$ Notation and terminology are
summarized in Section 2$.$ Supplementary results on supercyclicity and strong
stability required in the sequel are considered in Section 3$.$ The main
results are proved in Section 4, and synthesised in Corollary 4.3$.$ These
concern weak l-sequential supercyclicity and quasistability for power bounded
operators acting on a normed space$.$ Two applications are considered in
Section 5, characterising weakly l-sequentially supercyclic operators that are
not supercyclic (Corollary 5.2), and also characterising weakly l-sequentially
supercyclic isometries on Hilbert spaces (Corollary 5.5).

\section{Notation and Terminology }

Linear spaces in this paper are all over a generic scalar field $\FF$, which
is either $\RR$ or $\CC.$ Let $\R(L)$ and $\N(L)$ stand for range and kernel,
respectively, of a linear transformation $L$ between linear spaces$.$ Let
$\BX$ denote the normed algebra of all bounded linear transformations of a
normed space $\X$ into itself$.$ Elements of $\BX$ will be referred to as
operators$.$ Both norms, on $\X$ and the induced uniform norm on $\B[\X]$,
will be denoted by ${\|\cdot\|}.$ For each ${T\kern-1pt\in\kern-1pt\BX}$
take its power sequence $\{T^n\}_{n\ge0}.$ An operator $T$ is power bounded if
${\sup_n\|T^n\|<\infty}$, strongly stable if ${\|T^nx\|\to0}$ (i.e.,
${T^nx\to0}$) for every ${x\in\X}$ (notation: ${T^n\!\sconv O}$), and of class
$C_{1{\textstyle\cdot}\!}$ if ${\|T^nx\|\not\to0}$ (i.e., ${T^nx\not\to0}$)
for every nonzero ${x\in\X}.$ It is weakly stable if ${|f(T^nx)|\to0}$ for
\hbox{every} $f$ in the dual $\X^*$ of $\X$ (i.e., ${T^nx\wconv0}$) for every
${x\in\X}$ (notation: ${T^n\wconv O}).$ Strong stability implies weak
stability$.$ If $\X$ is a Banach space, then weak stability implies power
boundedness (by the Banach--Steinhaus Theorem)$.$ We say that an operator
${T\kern-1pt\in\kern-1pt\BX}$ on a normed space $\X$ is {\it weakly
quasistable}\/ if, for each ${x\in\X}$,
$$
{\liminf}_n|f(T^nx)|=0
\quad\hbox{for every}\quad
f\in\X^*\!.
$$
In particular, every weakly stable operator is \hbox{weakly quasistable}.

\vskip6pt
The {\it weak limit set}\/ $A^{-wl}$ of a subset $A$ of a normed space $\X$ is
the set of all weak limits of weakly convergent $A$-valued sequences, that is,
$$
A^{-wl}=\big\{x\in\X\!: x=w\hbox{\,-}\lim x_n\;\hbox{with}\;\,x_n\in A\big\}.
$$
\vskip-2pt\noi

\vskip6pt\noi
{\bf Remark 2.1.}
Let $A$ be a set in a normed space $\X.$ Its closure (in the norm topol\-ogy
of $\X$) is denoted by $A^-\!$ and its weak closure (
in the weak topology of
$\X$) is denoted by $A^{-w}\!$, so that $A$ is dense (in the norm topology)
or weakly dense (in the weak topology) if $A^-\!=\X$ or $A^{-w}\!=\X$,
respectively$.$ A set $A$ is weakly sequentially closed if every $A$-valued
weakly convergent sequence has its limit in $A.$ Let the weak sequential
closure of $A$ be denoted by $A^{-ws}\!$, which is the smallest (i.e., the
intersec\-tion of all) weakly sequentially closed subset of $\X$ including
$A.$ So $A$ is weakly sequentially dense if $A^{-ws}\!=\X.$ Consider the
above definition of the weak limit set $A^{-wl}\!$ of $A$ (the set of all
weak limits of weakly convergent $A$-valued sequences)$.$ We say that a set
$A$ is {\it weakly l-sequentially dense}\/ if $A^{-wl}\!=\X.$ It is known that
$$
A^-\!\sse A^{-wl}\!\sse A^{-ws}\!\sse A^{-w}\!,
$$ 
and the inclusions may be proper in general (see, e.g., \cite[pp.38,39]{Shk},
\cite[pp.259,260]{BM2}, \cite[pp.10,11]{Kub2})$.$ Therefore if $A$ is dense in
the norm topology (i.e., ${A^-\!=\X}$), then it is dense with respect to all
notions of denseness defined above$.$ Recall that if a set $A$ is convex, then
${A^-\!=A^{-w}}$ (see, e.g., \cite[Theorem 2.5.16]{Meg})$.$ Thus if $A$ is
convex, then the above inclusions become a chain of identities, and so for a
convex set the all the above notions of denseness coincide.

\vskip6pt
The orbit $\Oe_T(y)$ of a vector ${y\in\X}$ under an operator
${T\kern-1pt\in\kern-1pt\BX}$ is the set
$$
\Oe_T(y)={\bigcup}_{n\ge0}\{T^ny\}=\big\{T^ny\in\X\!:
\hbox{for every integer}\;n\ge0\big\}.
$$
A vector $y$ in $\X$ is a {\it cyclic vector}\/ for $T$ if the span of its
orbit is dense in $\X$:
$$
\big(\span\Oe_T(y)\big)^-\!=\X.
$$
An operator $T$ is a {\it cyclic operator}\/ if it has a cyclic vector$.$
Since ${\span\Oe_T(y)}$ is con\-vex,
${\big(\span\Oe_T(y)\big)^-\!=\big(\span\Oe_T(y)\big)^{-w}}\!$, and so the
notion of cyclicity is the same in the norm and in the weak topologies (i.e.,
{\it cyclicity coincides with weak cyclicity}\/).

\vskip6pt
Let $[x]=\span\{x\}$ denote the one-dimensional subspace of $\X$ spanned by a
sin\-gle\-ton $\{x\}$ at a vector ${x\in\X}.$ The projective orbit of a
nonzero vector ${y\in\X}$ under an operator ${T\in\BX}$ is the orbit
$\Oe_T([y])=\bigcup_{n\ge0}T^n([y])$ of the span of $\{y\}$:
$$
\Oe_T([y])=\big\{\alpha T^ny\in\X\!:\;
\hbox{for every}\;\alpha\in\FF\;\hbox{and every integer}\;n\ge0\big\}.
$$
A vector $y$ in $\X$ is a {\it supercyclic vector}\/ for $T$ if its projective
orbit $\Oe_T([y])$ is dense in $\X$; that is, if the orbit of the span of
$\{y\}$ is dense in $\X$:
$$
\Oe_T([y])^-\!=\X.
$$
Since the norm topology is metrizable, a nonzero vector $y$ in $\X$ is
supercyclic for an operator $T$ if and only if for every $x$ in $\X$ there
exists an $\FF$-valued sequence $\{\alpha_k\}_{k\ge0}$ (which depends on $x$
for each $y$ and consists of nonzero numbers if ${x\ne0}$) such that for some
sequence $\{T^{n_k}\}_{k\ge0}$ with entries from the sequence
$\{T^n\}_{n\ge0}$, the $\X$-valued sequence $\{\alpha_kT^{n_k}y\}_{k\ge0}$
converges to $x$ (in the norm topology): \hbox{for every ${x\in\X}$,}
$$
\alpha_kT^{n_k}y\to x
\qquad(\hbox{i.e.,}\quad
\|\alpha_kT^{n_k}y-x\|\to0).
$$
An operator $T$ is a {\it supercyclic operator}\/ if it has a supercyclic
vector.

\vskip6pt
The weak version of the above convergence criterion leads to the notion of
weak l-sequential supercyclicity$.$ A nonzero $y$ in $\X$ is a {\it weakly
l-sequentially supercyclic vector}\/ for $T$ if for every $x$ in $\X$ there is
an $\FF$-valued sequence $\{\alpha_k\}_{k\ge0}$ such that for some sequence
$\{T^{n_k}\}_{k\ge0}$ with entries from $\{T^n\}_{n\ge0}$, the $\X$-valued
sequence $\{\alpha_kT^{n_k}y\}_{k\ge0}$ converges weakly to $x.$ In other
words, if for every ${x\in\X}$,
$$
\alpha_kT^{n_k}y\wconv x
\qquad(\hbox{i.e.,}\quad
f(\alpha_kT^{n_k}y-x)\to0
\quad\hbox{for every}\quad
f\in\X^*).
$$
The $\FF$-valued sequence $\{\alpha_k\}_{k\ge0}$ depends on $x$ for each
weakly l-sequentially super\-cyclic vector $y$, and each
${\alpha_k=\alpha_k(x)}$ is nonzero whenever ${x\ne0}.$ The above displayed
convergence means that there are
$\{\alpha_k\}_{k\ge0}\!=\!\{\alpha_k(x)\}_{k\ge0}$ and
$\{T^{n_k}\}_{k\ge0}\!=\!\{T^{n_k(x)}\}_{k\ge0}$ (depending on $x$ for each
$y$ but not on $f$) such that, for each ${x\kern-1pt\in\kern-1pt\X}$,
$$
\alpha_k(x)f(T^{n_k(x)}y)\to f(x)
\quad\hbox{for every}\quad
f\in\X^*.
$$
Equivalently, the projective orbit $\Oe_T([y])$ of $y$ under $T\kern-1pt$ is
weakly l-sequentially dense in $\X$ in the sense that the weak limit set of
$\Oe_T([y])$ coincides with $\X$:
$$
\Oe_T([y])^{-wl}=\X.
$$
An operator $T$ is a {\it weakly l-sequentially supercyclic operator}\/ if it
has a weakly l-sequentially supercyclic vector.

\vskip6pt
The above notions are related as follows.
\vskip6pt\noi
\centerline{
{\ehsc supercyclic}
$\limply$
{\ehsc weakly l-sequentially supercyclic}
$\limply$
{\ehsc cyclic}
}
\vskip6pt\noi
(every supercyclic vector for $T$ is a weakly l-sequentially supercyclic
vector for $T$, which is a cyclic vector for $T$), and the converses fail$.$
For a comparison between these and further notions of cyclicity (including
hypercyclicity, weak hypercyclicity, weak sequential supercyclicity, and weak
supercyclicity) see, for instance, \cite[pp.38,39]{Shk},
\cite[pp.259,260]{BM2}, \cite[pp.159,232]{GP}, \cite[pp.50,51,54]{KD1},
\cite[pp.372,373,374]{KD2}$.$ Weak l-sequential supercyclicity is the central
theme in this paper$.$ It has been considered in \cite{BCS} (also in
\cite{BM1} implicitly), and discussed in \cite{Shk}, \cite{KD1}, and
\cite{KD3}$.$

\vskip6pt\noi
{\bf Remark 2.2.}
(a)
Every nonzero vector is trivially supercyclic for every operator on a
one-dimensional space, and all forms of cyclicity imply separability for
any normed space $\X.$ So we assume here that
$$
\hbox{\it all normed spaces are separable and have dimension greater than
one\/.}
$$
\vskip-2pt

\vskip6pt\noi
(b)
Take ${0\ne x\in\X}.$ If ${x\in\Oe_T([y])}$, then $\{T^{n_k}\}_{k\ge0}$ can be
viewed as a constant (infinite) sequence or, equivalently, as a one-entry
(thus finite) subsequence of $\{T^n\}_{n\ge0}$ (i.e., ${x=\alpha_0 T^{n_0}y}$
for some ${\alpha_0\ne0}$ and some ${n_0\ge0}).$ In this case, both forms of
supercyclicity (i.e., ${\alpha_kT^{n_k}y\to x}$ and
${\alpha_kT^{n_k}y\wconv x}$) coincide, where convergence of
$\{T^{n_k}y\}_{k\ge0}\kern-1pt$ means eventually reached$.$ Thus, in general,
we use the expression
$$
\hbox{\it $\{T^{n_k}\}_{k\ge0}$ is sequence with entries from\/
$\{T^n\}_{n\ge0}$}.
$$
However, if ${x\not\in\Oe_T([y])}$, this means
{\it $\{T^{n_k}\}_{k\ge0}$ is a subsequence of\/ $\{T^n\}_{n\ge0}$}.

\vskip6pt\noi
(c)
Let $T$ be a weakly l-sequentially supercyclic operator, and let $y$ be any
weakly l-sequentially supercyclic vector for $T.$ Then for every $x$ there
exists $\{\alpha_k(x)\}_{k\ge0}$ for which
${f(x)=\lim_k\alpha_k(x)f(T^{n_k(x)}y)}$ for every $f$, for some sequence
$\{T^{n_k(x)}\}_{k\ge0}$ with entries from $\{T^n\}_{n\ge0}.$ Note that
$$
\hbox{\it\/ if $T\kern-1pt$ is power bounded, then\/ $\{\alpha_k(x)\}_{k\ge0}$
cannot be uniformly bounded}\/.
$$
That is, if ${\sup_n\|T^n\|<\infty}$ (equivalently, if
${\sup_n\|T^nx\|<\infty}$ for every $x$ in a Banach space by the
Banach-Steinhaus Theorem), then it is not true that $\{\alpha_k(x)\}_{k\ge0}$
is such that $\,{\sup_x\sup_k|\alpha_k(x)|<\infty}.$ (Indeed, if $T$ is power
bounded and if $\{\alpha_k(x)\}_{k\ge0}$ is bounded for each $x$, then set
$\beta=\sup_n\|T^n\|$ and $\gamma(x)=\sup_k|\alpha_k(x)|$ for each $x$, so
that $|f(x)|\le{\beta\,\gamma(x)\|f\|\kern1pt\|y\|}$ for every $f$, and this
leads to a contradiction if ${\sup_x\gamma(x)<\infty}$ whenever $|f(x)|$ is
large enough, which can always be attained for an arbitrary $f$ since linear
functionals are surjective).

\section{Supercyclicity and Strong Stability}

The following result was proved in \cite[Theorem 2.1]{AB}.

\vskip6pt\noi
{\bf Proposition 3.1.}
\cite{AB}
{\it A power bounded supercyclic operator on a normed space is not of class}\/
$C_{1{\textstyle\cdot}}.$

\vskip6pt
This means that {\it if\/ ${T\kern-1pt\in\kern-1pt\BX}$ is such that\/
${\sup_n\|T^n\|<\infty}$, and if there exists a vector\/ ${y\in\X}$ such
that\/ ${\alpha_k(x)T^{n_k}y\to x}$ for every\/ ${x\in\X}$, for some\/
$\FF$-valued sequence}\/ $\{\alpha_k(x)\}_{k\ge0}$ (i.e., some sequence of
numbers that are nonzero whenever $x$ is nonzero) {\it and some sequence\/
$\{T^{n_k(x)}\}_{k\ge0}$ with entries from}\/ $\{T^n\}_{n\ge0}$ (both
depending on $x$ for each $y$), {\it then there exists a nonzero\/ ${z\in\X}$
for which}\/ ${T^nz\to0}$.

\vskip6pt
Using the above proposition, more was proved in \cite[Theorem 2.2]{AB},
namely, {\it if\/ ${T\kern-1pt\in\kern-1pt\BX}$ is such that\/
${\sup_n\|T^n\|<\infty}$, and if there exists a vector\/ ${y\in\X}$ such
that\/ ${\alpha_k(x)T^{n_k}y\to x}$ for every\/ ${x\in\X}$, for some sequence
of numbers\/ $\{\alpha_k(x)\}_{k\ge0}$ and some sequence
$\{T^{n_k(x)}\}_{k\ge0}$ with entries from\/ $\{T^n\}_{n\ge0}$, then\/
${T^nz\to0}$ for \hbox{every}}\/ ${z\in\X}.$ In other words, the result in
\cite[Theorem 2.2]{AB} may be read \hbox{as follows}.

\vskip6pt\noi
{\bf Proposition 3.2.}
\cite{AB}
{\it A power bounded supercyclic operator on a normed space is strongly
stable}\/$.$

\vskip6pt
Since supercyclicity implies weak l-sequential supercyclicity, and strong
stability implies weak stability, Proposition 3.2 prompts the question$:$
{\it does weak l-sequential supercyclicity imply weak stability for a power
bounded operator}$\,?$ The relationships around this question are summarized
in the following diagram.
\vskip2pt\noi
$$
\matrix{
& \Oe_T([y])^-=\X      & \limply & T^n\sconv O     \phantom{\Big|}  \cr
& \big\Downarrow       &         & \big\Downarrow                   \cr
& \Oe_T([y])^{-wl}=\X  & \query  & T^n\wconv O     \phantom{\Big|}. \cr}
$$
\goodbreak\noi
In other words,
$$
\matrix{
& \hbox{\ehsc$T$ is supercyclic} &\limply& \hbox{\ehsc$T$ is strongly stable}
\phantom{\Big|}                                                           \cr
& \big\Downarrow                 &       & \big\Downarrow                 \cr
& \hbox{\ehsc$T$ is weakly l-sequentially supercyclic}
                                 &\query & \hbox{\ehsc$T$ is weakly stable}
\phantom{\Big|}.                                                          \cr}
$$
\vskip2pt\noi
The above displayed question (originally posed in \cite{KD1}) remains open in
general$.$ There are, however, affirmative answers for particular classes of
operators$.$ For instance, if an operator is compact, then weak l-sequential
supercyclicity coincides with supercyclicity \cite[Theorem 4.1]{KD2}, which
implies strong stability (Proposition 3.2), which in turn implies weak
stability --- actually, in this case we get uniform stability
\cite[Theorem 4.2]{KD2}. We prove in Corollary 4.3 below the following
particular version of the above question for arbitrary power bounded
operators.
\vskip6pt\noi
\centerline{
{\ehsc $T$ is weakly l-sequentially supercyclic}
\quad$\limply$\quad
{\ehsc $T$ is weakly quasistable}.
}

\section{Main Results}

{From} now on let ${Y_T\sse\X}$ denote the set of all weakly l-sequentially
supercyclic vectors $y$ for an operator $T$ on a normed space $\X$,
$$
Y_T=\big\{y\in\X\!:\,\Oe_T([y])^{-wl}=\X\big\}
=\big\{
y\in\X\!:\,\forall\;x,\;\exists\;\{\alpha_k\}\;\ni\;\alpha_kT^{n_k}y\wconv x
\big\},
$$
so that ${0\not\in Y_T}$, and $T$ is weakly l-sequentially supercyclic if and
only if ${Y_T\ne\void}.$ Observe that ${Y_T\cup\0}$ is a $T$-invariant cone$.$
(Indeed, if ${y\in Y_T}$, then ${Ty\ne0}$ and every nonzero vector in
$\Oe_T([y])$ lies in $Y_T$ \cite[Lemma 5.1]{Kub2}; in particular, $Ty$ lies
in $Y_T$, and ${\alpha\kern1pty\in Y_T}$ for every nonzero ${\alpha\in\FF}$
trivially.)

\vskip6pt\noi
{\bf Theorem 4.1.}
{\it Suppose an operator\/ $T$ on a normed space\/ $\X$ is weakly
l-sequentially supercyclic$.$ Take an arbitrary}\/ ${y\in Y_T}$.

\begin{description}
\item{$\kern-6pt$\rm(a)}
{\it If there exists\/ an ${f_0\in\X^*\!}$ such that\/
${\liminf_n|f_0(T^ny)|>0}$, then for every\/ ${x\in\X}\kern1pt$ there exists
a\/ $\,${\rm bounded}$\,$ sequence of scalars\/
$\kern1pt\{\alpha_j(x)\}_{j\ge0}$ such that\/
${\alpha_j(x)f(T^{n_j(x)}y)\to f(x)}$ for every}\/ ${f\in\X^*}\!.$

\vskip4pt
\item{$\kern-6pt$\rm(b)}
{\it $\,$If for some\/ ${x\in\X}$ there exists a bounded scalar sequence\/
$\{\alpha_j(x)\}_{j\ge0}$ such that\/ ${\alpha_j(x)f(T^{n_j(x)}y)\to f(x)}$
for every\/ ${f\in\X^*}\!$, then there is a\/ $\,$\hbox{{\rm convergent}}$\,$
scalar sequence\/ $\{\alpha_k(x)\}_{k\ge0}$ such that\/
${\alpha_k(x)f(T^{n_k(x)}y)\to f(x)}$ for every}\/ ${f\in\X^*}$.

\vskip4pt
\item{$\kern-6pt$\rm(c)}
{\it $\,$Suppose the operator\/ $T$ is power bounded$.$ If for some nonzero\/
${x\in\X}$ there exists a convergent sequence of scalars\/
$\{\alpha_k(x)\}_{k\ge0}$, say\/ ${\lim_k\alpha_k(x)=\alpha(x)\in\FF}$, for
which\/ $\lim_k{\alpha_k(x)f(T^{n_k(x)}y)=f(x)}$ for every\/ ${f\in\X^*}\!$,
then\/ ${\alpha(x)\ne0}$ and\/ ${\alpha(x)\lim_kf(T^{n_k(x)}y)}=f(x)$ for
every}\/ ${f\in\X^*}\!$.

\vskip4pt
\item{$\kern-6pt$\rm(d)}
{\it Suppose\/ there is a functional\/ ${\alpha\!:\kern-1pt\X\!\to\FF}$ such
that\/ ${\alpha(x)\ne0}$ whenever\/ ${x\ne0}$ and, for each ${x\in\X}$,\/
${\alpha(x)\lim_k\kern-1ptf(T^{n_k(x)}y)=f(x)}$\/ for every\/ ${f\in\X^*}\!.$
Then ${\liminf_n|f(T^ny)|=0}$ for every}\/ ${f\in\X^*}\!$.
\end{description}

\proof
Let $T$ be a weakly l-sequentially supercyclic operator on a normed space
$\X$, and fix an arbritrary vector ${y\in Y_T\ne\void}$.

\vskip6pt\noi
(a)
Take an arbitrary $x$ in $\X.$ By weak l-sequential supercyclicity there is a
sequence $\{\alpha_j(x)\}_{j\ge0}$ such that
${\alpha_j(x)f(T^{n_j(x)}y)\to f(x)}$, and so
${|\alpha_j(x)|\,|f(T^{n_j(x)}y)|\to|f(x)|}$, for every
${f\kern-1pt\in\kern-1pt\X^*\!}.$ Suppose there is an
${f_0\kern-1pt\in\kern-1pt\X^*}$ for which
${\liminf_n\kern-1pt|f_0(T^ny)|\kern-1pt>\kern-1pt0}.\kern-1pt$ Then for every
$\veps$ in ${(0,\liminf_n\kern-1pt|f_0(T^ny)|)}$ there is an integer
${n_\veps>0}$ such that ${n\ge n_\veps}$ im\-plies ${\veps<|f_0(T^ny)|}.$
Thus, since
${\lim_j\kern-1pt|\alpha_j(x)|\kern1pt|f_0(T^{n_j(x)}y)|}
\kern-1pt=\kern-1pt|f_0(x)|$
is a real number, ${\limsup_j|\alpha_j(x)|\kern-1pt<\!\infty}.$ (Indeed,
${\veps\kern1pt\limsup_j\kern-1pt|\alpha_j(x)|}
<{\limsup_j\kern-1pt|\alpha_j(x)|\kern1pt|f_0(T^{n_j(x)}y)|}={|f_0(x)|}$
for $n_j(x)$ large enough --- recall$:$ $\{\alpha_j(x)\}_{j\ge0}$ depends on
$x$ but not on \hbox{$f$)$.$ So}
$$
\{\alpha_j(x)\}_{j\ge0}
\;\;\hbox{is bounded}.
$$
\vskip-2pt

\vskip6pt\noi
(b)
If $\{\alpha_j(x)\}_{j\ge0}$ is bounded, then (Bolzano--Weierstrass) there
exists a subsequence $\{\alpha_k(x)\}_{k\ge0}=\{\alpha_{j_k}(x)\}_{k\ge0}$ of
the sequence of scalars $\{\alpha_j(x)\}_{j\ge0}$ such that
$$
\{\alpha_k(x)\}_{k\ge0}
\;\;\hbox{converges}.
$$
Moreover, ${\alpha_k(x)f(T^{n_k(x)}y)\to f(x)}.$ Indeed, since
${\alpha_j(x)f(T^{n_j(x)}y)\to f(x)}$ for \hbox{every} $f$ in $\X^*\!$, the
convergence holds for every subsequence of $\{\alpha_j(x)f(T^{n_j(x)}y)\}$; in
particular, it holds for $\{\alpha_k(x)f(T^{n_k(x)}y)\}$, where
$\{T^{n_k(x)}\}_{k\ge0}=\{T^{n_{j_k}(x)}\}_{k\ge0}$ is the associated
subsequence of $\{T^{n_j}(x)\}_{j\ge0}$.

\vskip6pt\noi
(c)
Take an arbitrary ${x\in\X}.$ Let ${\alpha(x)\in\FF}$ be the limit of a
convergent sequence $\{\alpha_k(x)\}_{k\ge0}$ of scalars for which
${\alpha_k(x)f(T^{n_k(x)}y)\to f(x)}$ for every
${f\kern-1pt\in\kern-1pt\X^*}\!$ (where $\{\alpha_k(x)\}_{k\ge0}$ and
${n_k(x)}$ do not depend on $f).$ From the Hahn--Banach Theorem, if ${x\ne0}$,
then there exists ${f_1\kern-1pt\in\kern-1pt\X^*}$ such that ${f_1(x)\ne0}$
(as $\|x\|=$ $\sup_{\|f\|=1}|f(x)|\,).$ Since $T$ is power bounded, set
${\beta=\sup_n\|T^n\|}.$ Thus if \hbox{${x\ne0}$, then}
\begin{eqnarray*}
0<|f_1(x)|
&\kern-6pt=\kern-6pt&
{{\lim}_k|\alpha_k(x)|\kern1pt|f_1(T^{n_k(x)}y)|}                       \\
&\kern-6pt\le\kern-6pt&
{|\alpha(x)|\kern1pt{\limsup}_k|f_1(T^{n_k(x)}y)|}
\le{\beta\kern1pt|\alpha(x)|\kern1pt\|f_1\|\kern1pt\|y\|
<\infty},
\end{eqnarray*}
and hence ${\alpha(x)\ne0}.$ Summing up:
$$
x\ne0
\quad\hbox{implies}\quad
\alpha(x)\ne0.
$$
By definition, for ${x\ne0}$ the sequence $\{\alpha_k(x)\}_{k\ge0}$ is
such that ${\alpha_k(x)\ne0}$ for every ${k\ge0}.$ Thus if ${x\ne0}$, then
${0\ne\alpha_k(x)\to\alpha(x)\ne0}$, and so
${\alpha_k(x)^{-1}\to\alpha(x)^{-1}}.$ Since
${\alpha_k(x)f(T^{n_k(x)}y)\to f(x)}$, we get
${f(T^{n_k(x)}y)}
={\frac{\alpha_k(x)f(T^{n_k(x)}y)}{\alpha_k(x)}\to\frac{f(x)}{\alpha(x)}}$,
and so
$$
\alpha(x)\,{\lim}_kf(T^{n_k(x)}y)=f(x)
\quad\hbox{for every}\quad
f\in\X^*.
$$
\vskip-2pt

\vskip6pt\noi
(d)
We know from linear algebra that there is no injective linear functional on a
linear space of dimension greater than one (recall: we are assuming
${\dim\X\kern-1pt>1}).$ Take an arbitrary ${0\ne f\in\X^*}\!.$ Now take a
nonzero ${x_0\in\kern-1pt\X}$ in the kernel \hbox{$\N(f)\ne\0$ of $f$}.

\vskip6pt\noi
(d$_1$)
If ${x_0\in\Oe_T([y])}$, then $x_0={\alpha_0T^{n_0}y}$ for some scalar
${\alpha_0\ne0}$ and some nonnegative integer $n_0.$ So
${\alpha_0f(T^{n_0}y)=f(x_0)=0}.$ Then ${f(T^{n_0}y)=0}.$ But this cannot hold
for every ${f\in\X^*}$ (since there is no ${0\ne z\in\X}$ for which ${f(z)=0}$
for every ${f\in\X^*}$.) Thus for each nonzero ${f\in\X^*}$ there exists a
nonzero ${x_0\in\N(f)\backslash\Oe([y])}$.

\vskip6pt\noi
(d$_2$)
Then suppose ${x_0\not\in\Oe([y])}.$ Therefore
${\alpha(x_0)\lim_kf(T^{n_k(x_0)}y)=f(x_0)}$ for some subsequence
$\{T^{n_k(x_0)}\}_{k\ge0}$ of $\{T^n\}_{n\ge0}.$ Since ${\alpha(x_0)\ne0}$,
and since
$$
\alpha(x_0)\,{\lim}_kf(T^{n_k(x_0)}y)=f(x_0)=0,
$$
we get ${\lim_kf(T^{n_k(x_0)}y)=0}.$ Then for an arbitrary ${0\ne f\in\X^*}\!$
there is a subsequence of $\{f(T^ny)\}$ converging to zero$.$ So
${\liminf_n|f(T^ny)|\kern-1pt=\kern-1pt0}$ for every
${f\kern-1pt\in\kern-1pt\X^*}\!.\!$                                      \qed

\vskip6pt
Before harvesting the corollaries we need the following lemma.

\vskip6pt\noi
{\bf Lemma 4.2.}
{\it Let\/ $T$ be a power bounded weakly l-sequentially supercyclic operator
 on a normed space}\/ $\X.$
\begin{description}
\item{$\kern-4pt$\rm (a)}
{\it $T$ is weakly stable if and only if\/ ${T^ny\wconv0}$ for every\/
${y\in Y_T}$
\vskip0pt\noi
$($i.e., if and only if\/ ${\lim_n|f(T^ny)|=0}$ for every}\/ ${f\in\X^*}\!$,
for every\/ ${y\in Y_T}$).
\vskip4pt
\item{$\kern-4pt$\rm (b)}
{\it $T$ is weakly quasistable if and only if\/ ${\liminf_n|f(T^ny)|=0}$ for
every\/ ${f\in\X^*}\!$, for every}\/ ${y\in Y_T}$.
\end{description}

\proof
That the set of weakly supercyclic vectors is dense in the norm topology was
shown in \cite[Proposition 2.1 ]{San}$.$ This has been extended to the set of
weakly l-sequentially supercyclic vectors in \cite[Theorem 5.1]{Kub2}, so
that
$$
Y_T\ne\void
\quad\limply\quad
Y_T^-\!=\X
\quad
\hbox{(and so $Y_T^{-wl}\!=\X$)}.
$$
Take an arbitrary $x$ in $\X.$ Since ${Y_T^-=\X}$, there exists a $Y_T$-valued
sequence $\{y_k\}$ such that ${\|y_k\!-x\|\to0}.$ If ${\lim_n|f(T^ny)|=0}$ (or
${\liminf_n|f(T^ny)|=0}$) for every $f$ in $\X^*\!$ and for every $y$ in $Y_T$,
then ${\lim_n|f(T^ny_k)|=0}$ (or ${\liminf_n|f(T^ny_k)|=0}$) for every $f$ in
$\X^*\!$ and for each integer $k.$ Thus, since $T$ is power bounded, and since
$$
|f(T^nx)|\le|f(T^n(y_k-x))|+|f(T^ny_k)|
\le\|f\|\,{\sup}_n\|T^n\|\kern1pt\|y_k-x\|+|f(T^ny_k)|
$$
for every ${f\in\X^*\!}$ and ${x\in\X}$, we get ${\lim_n|f(T^nx)|\to0}$ (or
${\liminf_n|f(T^nx)|=0}$) for every ${f\in\X^*}$ and every ${x\in\X}.$ Hence
${\lim_n|f(T^ny)|=0}$ (or ${\liminf_n|f(T^ny)|=0}$) for every ${f\in\X^*}$ and
every ${y\in Y_T}$ implies $T$ is weakly stable (or weakly quasistable), and
the converse (in both cases) is trivial.                                 \qed

\vskip6pt
The result below consists of an improvement over \cite[Theorem 6.2]{KD1},
ensuring that weak quasistability always happens for weakly l-sequentially
supercyclic power bounded operators acting on arbitrary normed spaces.

\vskip6pt\noi
{\bf Corollary 4.3.}
{\it Every weakly l-sequentially supercyclic power bounded operator on a
normed space is weakly quasistable}\/.

\proof
Let $T$ be a power bounded weakly l-sequentially supercyclic operator on a
normed space $\X.$ Take an arbitrary weakly l-sequentially supercyclic vector
$y$ in $Y_T.$ Theorem 4.1 ensures that if there exists ${f_0\in\X^*}$ such
that ${\liminf_n|f_0(T^ny)|>0}$, then ${\liminf_n|f(T^ny)|=0}$ for every
${f\in\X^*}\!$, which is a contradiction$.$ So if a power bounded operator
$T$ is weakly l-sequentially supercyclic, then ${\liminf_n|f(T^ny)|}=0$ for
every ${f\kern-1pt\in\kern-1pt\X^*}\!$, for every
${y\kern-1pt\in\kern-1ptY_T}.$ Thus $T$ is weakly quasistable by
\hbox{Lemma 4.2(b).\!\!\!\qed}

\vskip6pt\noi
{\bf Corollary 4.4.}
{\it If\/ $T$ is a power bounded operator on a normed space, then}\/
$$
\Oe_T([y])^{-wl}=\X
\;\;\hbox{\it for some}\;\;
y\in\X
\quad\limply\quad
0\in\Oe_T(x)^{-wl}
\;\;\hbox{\it for every}\;\;
x\in\X.
$$
\vskip-2pt\noi

\proof
Suppose ${Y_T\ne\void}$ and take an arbitrary ${y\in Y_T}.$ Now take an
arbitrary ${x\in\X}.$ If $T$ is a power bounded operator on $\X$, then
Corollary 4.3 says that ${\Oe_T([y])^{-wl}=\X}$ implies
${\liminf_n|f(T^nx)|=0}$ for every ${f\in\X^*}\!.$ This means that there is
an $\X$-valued subsequence $\{T^{n_j}x\}_{j\ge0}$ of $\{T^nx\}_{n\ge0}$ (which
actually is an $\Oe_T(x)$-valued sequence that does not depend on $f$) such
that ${\lim_j|f(T^{n_j}x)|=0}$ for every ${f\in\X^*}\!.$ That is,
${T^{n_j}x\wconv0}$, so that zero is a weak limit of the weakly convergent
$\Oe_T(x)$-valued sequence $\{T^{n_j}x\}_{j\ge0}$, which in turn means that
${0\in\Oe_T(x)^{-wl}}.$                                                \qed

\section{Applications}

A Radon--Riesz space is a normed space $\X$ for which the following property
holds: an $\X$-valued sequence $\{x_n\}_{n\ge0}$ converges strongly (i.e.,
converges in the norm topology) if and only if it converges weakly and the
sequence of norms $\{\|x_n\|\}_{n\ge0}$ converges to the norm of the limit
(see, e.g., \cite[Definition 2.5.26]{Meg})$.$ In oder words, a normed space
$\X$ is a Radon-Riesz space if, for an arbitrary $\X$-valued
\hbox{sequence $\{x_n\}$},
$$
x_n\to x
\quad\iff\quad
\big\{x_k\wconv x
\;\;\,\hbox{and}\,\;\;\|x_n\|\to\|x\|\big\}.
$$
Hilbert spaces are Radon--Riesz spaces (cf$.$
\hbox{\cite[Problem 20, p.13]{Hal}).}

\vskip6pt\noi
{\bf Lemma 5.1.}
{\it Let\/ $\X$ be a Radon--Riesz space$.$ If a\/n $\X$-valued sequence\/
$\{x_n\}_{n\ge0}$ converges weakly to\/ ${x\in\X}$ but not strongly\/
$($i.e., ${x_n\wconv x}$ and\/ ${x_n\not\to x})$, then}\/
$$
\|x\|\le{\liminf}_n\|x_n\|<{\limsup}_n\|x_n\|.
$$

\proof
If $\X$ is a normed space, then ${x_n\wconv x}$ implies
$\|x\|\le{\liminf}_n\|x_n\|$ (see, e.g., \cite[Proposition 46.1]{Heu} ---
compare with \cite[Problem 21, p.14]{Hal} --- norm is weakly lower
semicontinuous)$.$ Thus if ${x_n\not\to x}$ and $\X$ is a Radon--Riesz space,
then ${\|x_n\|\not\to\|x\|}$ by the above displayed equivalence, and this
implies that $\{\|x_n\|\}_{n\ge0}$ does not converge, and hence
${\liminf_n\|x_n\|}<{\limsup_n\|x_n\|}.$ This concludes the proof.     \qed

\vskip6pt
As Corollary 4.3 supplied an improvement over \cite[Theorem 6.2]{KD1}, the
result below supplies an improvement over \cite[Theorem 5.1]{KD1}, in the
sense that it offers a sharper condition for distinguishing weak l-sequential
supercyclicity from (strong) supercyclicity$.$ Suppose $T$ is a power bounded
weakly l-sequentially supercyclic op\-erator on a Radon--Riesz space $\X$
which is not supercyclic$.$ Take an arbitrary $y$ in $Y_T$ so that, for each
${x\in\X}$, there is a sequence of scalars $\{\alpha_k(x)\}_{k\ge0}$ for which
${\alpha_k(x)T^{n_k(x)}y\wconv x}$ and ${\alpha_k(x)T^{n_k(x)}y\not\to x}.$
Theorem 4.1(b,c,d) and Lemma 4.2(b) show that, if $\{\alpha_k(x)\}_{k\ge0}$ is
bounded for each ${x\in\X}$, then $T$ is weakly quasistable, and by
Theorem 4.1(a) we saw in Corollary 4.3 that weak quasistability for $T$ always
holds$.$ {\it Does boundedness\/ $($thus convergence$)$ of\/
$\{\alpha_k(x)\}_{k\ge0}$ \hbox{always holds as well}$\,?$}

\vskip6pt\noi
{\bf Corollary 5.2.}
{\it Let\/ $T$ be a power bounded operator on a Radon--Riesz space$.$ Suppose\/
$T$ is weakly l-sequentially supercyclic but not supercyclic$.$ Take any
vector\/ $y$ in\/ $Y_T$ so that for every\/ $x$ in ${\X\backslash\Oe_T([y])}$
there exists a sequence of scalars\/ $\{\alpha_k(x)\}_{k\ge0}$ for which\/
${\alpha_k(x)\,T^{n_k(x)}y\wconv x}$ for some subsequence\/
$\{T^{n_k(x)}\}_{k\ge0}$ of\/ $\{T^n\}_{n\ge0}.$ If\/ $\{\alpha_k(x)\}_{k\ge0}$
is bounded for each\/ ${x\in\X\backslash\Oe_T([y])}$, then for every nonzero\/
${f\in\X^*}$
$$
{\limsup}_j|f(T^{n_j}y)|<\|f\|\kern1pt{\limsup}_j\|T^{n_j}y\|        \eqno(*)
$$
for some subsequence\/ $\{T^{n_j}\}_{j\ge0}$ of}\/ $\{T^n\}_{n\ge0}$.

\proof
Let $T$ be an operator on a normed space $\X.$ If it is weakly l-sequentially
supercyclic but not supercyclic, then take any weakly l-sequentially
supercyclic vector $y$ in $Y_T$, and recall that it is not a supercyclic
vector (there is no supercyclic vector)$.$ Take an arbitrary $x$ in
${\X\backslash\Oe_T([y])}.$ Thus there is a sequence of scalars
$\{\alpha_j(x)\}_{j\ge0}$ such that ${\alpha_j(x)\,T^{n_j(x)}y\wconv x}$ for
some subsequence $\{T^{n_j(x)}\}_{j\ge0}$ of $\{T^n\}_{n\ge0}.$ Since
$\{\alpha_j(x)\}_{j\ge0}$ is bounded for every ${x\in\X}$ and $T$ is power
bounded, Theorem 4.1(b,c) says that, in this case, there is a nonzero constant
scalar sequence, say $\alpha_k(x)=\alpha(x)$ for every ${k\ge0}$, for which
${\alpha(x)\,T^{n_k(x)}y\wconv x}.$ Hence, for each
${x\in\X\backslash\Oe_T([y])}$,
$$
|f(x)|=|\alpha(x)|\,{\lim}_k|f(T^{n_k(x)}y)|
\quad\hbox{for every}\quad
f\in\X^*.                                                           \eqno(1)
$$
As $y$ is not supercyclic, there is an ${x_0\ne0}$ in
${\X\backslash\Oe_T([y])}$ such that ${\alpha_i\,T^{n_i}y\not\to x_0}$ for
every sequence of scalars $\{\alpha_i\}_{i\ge0}$ and every subsequence
$\{T^{n_i}\}_{i\ge0}$ of \hbox{$\{T^n\}_{n\ge0}.$ So}
$$
\alpha(x_0)\,T^{n_k(x_0)}y\wconv x_0
\quad\;\hbox{and}\;\quad
\alpha(x_0)\,T^{n_k(x_0)}y\not\to x_0.
$$
As $\X$ is a Radon--Riesz space, it follows by Lemma 5.1 that
$$
\|x_0\|
\le{\liminf}_k\|\alpha(x_0)\,T^{n_k(x_0)}y\|
<{\limsup}_k\|\alpha(x_0)\,T^{n_k(x_0)}y\|.                         \eqno(2)
$$
Now suppose $(*)$ fails$.$ That is, suppose there exists a nonzero
${f_1\in\X^*}$ such that
$$
{\limsup}_i|f_1(T^{n_i}y)|=\|f_1\|{\limsup}_i\|T^{n_i}y\|           \eqno(3)
$$
for every subsequence $\{T^{n_i}\}_{i\ge0}$ of $\{T^n\}_{n\ge0}.$
Then, by (1), (2), and (3),
\begin{eqnarray*}
|\alpha(x_0)|\,{\limsup}_k|f_1(T^{n_k(x_0)}y)|
&\kern-6pt=\kern-6pt&
|\alpha(x_0)|\,{\lim}_k|f_1(T^{n_k(x_0)}y)|                              \\
&\kern-6pt=\kern-6pt&
|f_1(x_0))|\le\|f_1\|\,\|x_0\|                                           \\
&\kern-6pt<\kern-6pt&
\|f_1\|\,{\limsup}_k\|\alpha(x_0)T^{n_k(x_0)}y\|                         \\
&\kern-6pt=\kern-6pt&
|\alpha(x_0)|\,{\limsup}_k|f_1(T^{n_k(x_0)}y|,
\end{eqnarray*}
which is a contradiction$.$ Thus $(*)$ holds.                          \qed

\vskip6pt
To proceed we need the following properties of weak l-sequential
supercyclicity.

\vskip6pt\noi
{\bf Lemma 5.3.}
{\it A weakly l-sequentially supercyclic operator has dense range}\/.

\proof
Take an operator ${T\kern-1pt\in\kern-1pt\BX}$ on a normed space $\X.$ If a
vector $y$ in $\X$ is a weakly l-sequentially supercyclic vector for $T$, then
any vector in the projective orbit $\Oe_T(\span\{y\})$ is again weakly
l-sequentially supercyclic \cite[Lemma 5.1]{Kub2}, and so is the vector $Ty.$
Thus the projective orbit of $Ty$ is included in the range of $T$:
$$
\Oe_T(\span\{Ty\})
=\big\{\alpha T^ny\in\X\!:\,\alpha\in\FF,\,n\ge1\big\}\sse\R(T)\sse\X.
$$
As $Ty$ is a weakly l-sequentially supercyclic vector for $T$, then
$\Oe_T(\span\{Ty\})$ is weakly l-sequentially dense in $\X$, and so is $\R(T)$
(i.e., ${\R(T)^{-wl}\!=\X}$)$.$ Since $\R(T)$ is a linear manifold of $\X$, it
is convex, and therefore it is dense (in the norm topology) in $\X$ (i.e.,
${\R(T)^-\!=\X}$) according to Remark 2.1.                               \qed

\vskip6pt\noi
{\bf Lemma 5.4.}
{\it A weakly l-sequentially supercyclic isometry on a Banach space is
surjective\/$.\!$ So a weakly l-sequentially supercyclic isometry on a Hilbert
space is \hbox{unitary}.}

\proof
Isometries on Banach spaces have closed range (because they are
bounded and bounded below)$.$ Thus apply Lemma 5.3 and get a surjective
isometry$.$ Surjective isometries on Hilbert spaces are precisely the unitary
operators. \qed

\vskip6pt
It is known that there are weakly l-sequentially supercyclic unitary operators
\cite[Example 3.6]{BM1}, \cite[pp.10,12]{BM1},
\cite[Proposition 1.1, Theorem 1.2]{Shk})$.$ A unitary operator (acting on a
Hilbert space) is singular-continuous if its scalar spectral measure is
singular-continuous with respect to normalized Lebesgue measure on the
$\sigma$-algebra of Borel subsets of the unit circle.

\vskip6pt\noi
{\bf Corollary 5.5.}
{\it If an isometry on a Hilbert space is weakly l-sequentially supercyclic,
then it is a weakly quasistable singular-continuous unitary operator}\/.

\proof
Take a weakly l-sequentially supercyclic isometry on a Hilbert space$.$
Isometries are power bounded, so the conclusion of weak quasistability
follows from Cor\-ollary 4.3$.$ Lemma 5.4 says that the isometry must be
unitary$.$ But every weakly l-se\-quentially supercyclic unitary operator is 
singular-continuous \hbox{\cite[Theorem 4.2]{Kub1}.}                     \qed

\vskip6pt
It is also known that there are singular-continuous unitary operators that are
weakly stable, and singular-continuous unitary operators that are not weakly
stable \cite[Propositions 3.2 and 3.3]{Kub1}$.$ We close the paper with a
question$.$ 
$$
\hbox{\it Is there a singular-continuous unitary operator that is not weakly
quasistable$\,?$}
$$

\vskip10pt\noi
\bibliographystyle{amsplain}

\begin{thebibliography}{10}

\bibitem{AB}
S.I. Ansari and P.S. Bourdon,
{\it Some properties of cyclic operators}\/,
Acta Sci. Math. (Szeged) {\bf 63} (1997), 195--207.

\bibitem{BM1}
F. Bayart and E. Matheron,
{\it Hyponormal operators, weighted shifts and weak forms of
supercyclicity}\/,
Proc. Edinb. Math. Soc. {\bf 49} (2006), 1--15.

\bibitem{BM2}
F. Bayart and E. Matheron,
{\it Dynamics of Linear Operators}\/,
Cambridge University Press, Cambridge, 2009.

\bibitem{BCS}
J. Bes, K.C. Chan, and R. Sanders,
{\it Every weakly sequentially hypercyclic shift is norm hypercyclic}\/.
Math. Proc. R. Ir. Acad. {\bf 105A} (2005), 79--85.

\bibitem{GP}
K.-G. Grosse-Erdmann and A. Peris,
{\it Linear Chaos}\/,
Springer, London, 2011.

\bibitem{Hal}
P.R. Halmos,
{\it A Hilbert Space Problem Book}\/,
2nd edn. Springer, New York, 1982.

\bibitem{Heu}
H.G. Heuser,
{\it Functional Analysis}\/,
Wiley, Chichester, 1982.

\bibitem{Kub1}
C.S. Kubrusly,
{\it Singular-continuous unitaries and weak dynamics}\/,
Math. Proc. R. Ir. Acad. {\bf 116A} (2016), 45--56.

\bibitem{Kub2}
C.S. Kubrusly,
{\it Denseness of sets of supercyclic vectors}\/,
Math. Proc. R. Ir. Acad. {\bf 120A} (2020), 7--18.

\bibitem{KD1}
C.S. Kubrusly and B.P. Duggal,
{\it On weak supercyclicity I}\/,
Math. Proc. R. Ir. Acad. {\bf 118A} (2018), 47--63.

\bibitem{KD2}
C.S. Kubrusly and B.P. Duggal,
{\it On weak supercyclicity II}\/,
Czechoslovak Math. J. {\bf 68(143)} (2018), 371--386.

\bibitem{KD3}
C.S. Kubrusly and B.P. Duggal,
{\it Weakly supercyclic power bounded operators of class}\/ $C_{1.}$,
Adv. Math. Sci. Appl. {\bf 30} (2021), 571--585.

\bibitem{Meg}
R. Megginson,
{\it An Introduction to Banach Space Theory}\/,
Springer, New York, 1998.

\bibitem{San}
R. Sanders,
{\it Weakly supercyclic operators}\/,
J. Math. Anal. Appl. {\bf 292} (2004), 148--159.

\bibitem{Shk}
S. Shkarin,
{\it Non-sequential weak supercyclicity and hypercyclicity}\/,
J. Funct. Anal. {\bf 242} (2007), 37--77.

\end{thebibliography}

\end{document}